\documentclass[12pt]{amsproc}

\author{Emerson de Melo}
\address{Department of Mathematics, University of Bras\'ilia, Bras\'ilia-DF 70910-900, Brazil}
\email{emerson@mat.unb.br}
\author{Pavel Shumyatsky}
\address{Department of Mathematics, University of Bras\'ilia,
Bras\'ilia-DF, 70910-900, Brazil}
\email{pavel@unb.br}

\thanks{The first author was supported by FEMAT; The second author was supported by FAPDF and CNPq-Brazil}

\title{Fitting subgroup and nilpotent residual of fixed points}
\date{2017}

\newtheorem{theorem}{\sc Theorem}[section]
\newtheorem{lemma}[theorem]{\sc Lemma}

\begin{document}

\begin{abstract}
Let $q$ be a prime and $A$ an elementary abelian group of order at least $q^3$ acting by automorphisms on a finite $q'$-group $G$. It is proved that if $|\gamma_{\infty}(C_{G}(a))|\leq m$ for any $a\in A^{\#}$, then the order of  $\gamma_{\infty}(G)$ is $m$-bounded. If $F(C_{G}(a))$ has index at most $m$ in $C_G(a)$ for any $a \in A^{\#}$, then the index of $F_2(G)$ is $m$-bounded.
\end{abstract}

\maketitle

\section{Introduction}

Suppose that a finite group $A$ acts by automorphisms on a finite group $G$.  The action is coprime if the groups $A$ and $G$ have coprime orders. We denote by $C_G(A)$ the set $$\{g\in G\ |\ g^a=g \ \textrm{for all} \ a\in A\},$$ the centralizer of $A$ in $G$ (the fixed-point subgroup). In what follows we denote by $A^\#$ the set of nontrivial elements of $A$. It has been known that centralizers of coprime automorphisms have strong influence on the structure of $G$.

Ward showed that if $A$ is an elementary abelian $q$-group of rank at least 3 and if $C_G(a)$ is nilpotent for any $a\in A^\#$, then the group $G$ is nilpotent \cite{War}. Later the second author showed that if, under these hypotheses, $C_G(a)$ is nilpotent of class at most $c$ for any $a\in A^\#$, then the group $G$ is nilpotent with $(c,q)$-bounded nilpotency class \cite{Sh1}. Throughout the paper we use the expression ``$(a,b,\dots )$-bounded'' to abbreviate ``bounded from above in terms of  $a,b,\dots$ only''. Subsequently the above result was extended to the case where $A$ is not necessarily abelian. Namely, it was shown in \cite{Eme1} that if $A$ is a finite group of prime exponent $q$ and order at least $q^3$ acting on a finite $q'$-group $G$ in such a manner that $C_G(a)$ is nilpotent of class at most $c$ for any $a\in A^{\#}$, then $G$ is nilpotent with class bounded solely in terms of $c$ and $q$. Many other results illustrating the influence of centralizers of automorphisms on the structure of $G$ can be found in \cite{khukhro}.

In the present article we address the case where $A$ is an elementary abelian $q$-group of rank at least 3 and $C_G(a)$ is ``almost" nilpotent for any $a\in A^\#$. Recall that the nilpotent residual of a finite group $G$ is the intersection of all terms of the lower central series of $G$. This will be denoted by $\gamma_\infty(G)$. One of the results obtained in \cite{Eme2} says that if $A$ and $G$ are as above and $\gamma_\infty(C_G(a))$ has order at most $m$ for any $a\in A^\#$, then the order of $\gamma_\infty(G)$ is $(m,q)$-bounded. The purpose of the present article is to obtain a better result by showing that the order of $\gamma_\infty(G)$ is $m$-bounded and, in particular, the order of $\gamma_\infty(G)$ can be bounded by a number independent of the order of $A$.

\begin{theorem}\label{main1}
Let $q$ be a prime and $m$ a positive integer. Let $A$ be an elementary abelian group of order at least $q^3$ acting by automorphisms on a finite $q'$-group $G$. Assume that $|\gamma_{\infty}(C_{G}(a))|\leq m$ for any $a\in A^{\#}$. Then $|\gamma_{\infty}(G)|$ is $m$-bounded.   
\end{theorem}

Further, suppose that the Fitting subgroup of $C_G(a)$ has index at most $m$ in $C_G(a)$ for any $a\in A^\#$. It was shown in \cite{Sh} that under this assumption the index of the Fitting subgroup of $G$ is $(m,q)$-bounded. In view of Theorem \ref{main1} it is natural to conjecture that in fact the index of the Fitting subgroup of $G$ can be bounded in terms of $m$ alone. We have not been able to confirm this. Our next result should be regarded as an evidence in favor of the conjecture. Recall that the second Fitting subgroup $F_2(G)$ of a finite group $G$ is defined as the inverse image of $F(G/F(G))$, that is, $F_2(G)/F(G)=F(G/F(G))$. Here $F(G)$ stands for the Fitting subgroup of $G$.

\begin{theorem}\label{main2}
Let $q$ be a prime and $m$ a positive integer. Let $A$ be an elementary abelian group of order at least $q^3$ acting by automorphisms on a finite $q'$-group $G$. Assume that $F(C_{G}(a))$ has index at most $m$ in $C_G(a)$ for any $a \in A^{\#}$. Then the index of $F_2(G)$ is $m$-bounded.   
\end{theorem}

In the next section we give some helpful lemmas that will be used in the proofs of the above results. Section 3 deals with the proof of Theorem \ref{main2}. In Section 4 we prove Theorem \ref{main1}.
\section{Preliminaries}
If $A$ is a group of automorphisms of a group $G$, the subgroup generated by elements of the form $g^{-1}g^\alpha$ with $g\in G$ and $\alpha\in A$ is denoted by $[G,A]$. The subgroup $[G,A]$ is an $A$-invariant normal subgroup in $G$. Our first lemma is a collection of well-known facts on coprime actions (see for example \cite{GO}). Throughout the paper we will use it without explicit references.
\begin{lemma}\label{111} Let $A$ be a group of automorphisms of a finite group $G$ such that $(|G|,|A|)=1$. Then
\begin{enumerate}
\item[i)] $G=[G,A]C_{G}(A)$.
\item[ii)] $[G,A,A]=[G,A]$.
\item[iii)] $A$ leaves invariant some Sylow $p$-subgroup of $G$ for each prime $p\in\pi(G)$.
\item[iv)] $C_{G/N}(A)=C_G(A)N/N$ for any $A$-invariant normal subgroup $N$ of $G$.
\item[v)] If $A$ is a noncyclic elementary abelian group and $A_1,\dots,A_s$ are the maximal subgroups in $A$, then $G=\langle C_G(A_1),\ldots,C_G(A_s)\rangle$. Furthermore, if $G$ is nilpotent, then $G=\prod_i C_G(A_i)$.
\end{enumerate}
\end{lemma}

The following lemma was proved in \cite{almo}. The case where the group $G$ is soluble was established in Goldschmidt \cite[Lemma 2.1]{gold}. 
\begin{lemma}\label{golds} Let $G$ be a finite group acted on by a finite group $A$ such that $(|A|,|G|)=1$. Then $[G,A]$ is generated by all nilpotent subgroups $T$ such that $T=[T,A]$.
\end{lemma}

\begin{lemma}\label{L1}
Let $q$ be a prime and $A$ an elementary abelian group of order at least $q^2$ acting by automorphisms on a finite $q'$-group $G$. Let $A_1,\dots,A_s$ be the subgroups of index $q$ in $A$. Then $[G,A]$ is generated by the subgroups $[C_G(A_i),A]$.
\end{lemma}
\begin{proof} If $G$ is abelian, the result is immediate from Lemma \ref{111}(v) since the subgroups $C_G(A_i)$ are $A$-invariant. If $G$ is nilpotent, the result can be obtained by considering the action of $A$ on the abelian group $G/\Phi(G)$. Finally, the general case follows from the nilpotent case and Lemma \ref{golds}.
\end{proof}

The following lemma is an application of the three subgroup lemma.
\begin{lemma}\label{lz}
Let  $A$ be a group of automorphisms of a finite group $G$ and let $N$ be a normal subgroup of $G$ contained in $C_G(A)$. Then $[[G,A],N]=1$. In particular, if $G=[G,A]$, then $N\leq Z(G)$.
\end{lemma}
\begin{proof} Indeed, by the hypotheses, $[N,G,A]=[A,N,G]=1$. Thus, $[G,A,N]=1$ and the lemma follows.
\end{proof}

In the next lemma we will employ the fact that if $A$ is any coprime group of automorphisms of a finite simple group, then $A$ is cyclic (see for example \cite{GS}). We denote by $R(H)$ the soluble radical of a finite group $H$, that is, the largest normal soluble subgroup of $H$.

\begin{theorem}\label{solu}
Let $q$ be a prime and $m$ a positive integer such that $m<q$. Let $A$ be an elementary abelian group of order $q^2$ acting on a finite $q'$-group $G$ in such a way that the index of $R(C_G(a))$ in $C_G(a)$ is at most $m$ for any $a\in A^{\#}$. Then $[G,A]$ is soluble.
\end{theorem}
\begin{proof}
We argue by contradiction. Choose a counter-example $G$ of minimal order. Then $G=[G,A]$ and $R(G)=1$. Suppose that $G$ contains a proper normal $A$-invariant subgroup $N$. Since $[N,A]$ is subnormal, we conclude that $[N,A]=1$ and so $N=C_N(A)$. In that case by Lemma \ref{lz} $N$ is central and, in view of $R(G)=1$, we have a contradiction. 

Hence $G$ has no proper normal $A$-invariant subgroups and so $G=S_1\times\dots S_l$, where $S_i$ are isomorphic nonabelian simple subgroups transitively permuted by $A$. We will prove that under these assumptions $G$ has order at most $m$.

If $l=1$, then $G$ is a simple group and so $G=C_G(a)$ for some $a\in A^\#$. In this case we conclude that $G$ has order at most $m$ by the hypotheses. Suppose that $l\neq1$ and so $l=q$, or $l=q^2$.

In the first case $G=S\times S^a\times\dots\times S^{a^{q-1}}$ for some  $a\in A$ and there exists $b\in A$ such that $S^b=S$. Here $S=S_1$. We see that $C_G(a)$ is the ``diagonal" of the direct product. In particular $C_G(a)\cong S$ is a simple group and so $C_G(a)$ is of order at most $m$. Since $m<q$ and $b$ leaves $C_G(a)$ invariant we conclude that $C_G(a)\leq C_G(b)$. Combining this with the fact that $b$ stabilizes all simple factors, we deduce that $b$ acts trivially on $G$. It follows that $|G|\leq m$.

Finally suppose that $G$ is a product of $q^2$ simple factors which are transitively permuted by $A$. For each $a\in A$ we see that $C_G(a)$ is a product of $q$ ``diagonal" subgroups. In particular,  $C_G(a)$ contains a direct product of $q$ nonabelian simple groups. This is a contradiction since $[C_G(a):R(C_G(a))]$ is at most $m$ and $m<q$. 

This proves that $G$ has order at most $m$. Then of course $A$ acts trivially on $G$. We conclude that $[G,A]=1$ that is a contraction and the proof is complete.
\end{proof}

\section{Proof of Theorem  \ref{main2}}

Assume the hypothesis of Theorem \ref{main2}. Thus, $A$ is an elementary abelian group of order at least $q^3$ acting on a finite $q'$-group $G$ in such a manner that $F(C_{G}(a))$ has index at most $m$ in $C_G(a)$ for any $a\in A^{\#}$. We wish to show that $F_2(G)$ has $m$-bounded index in $G$. It is clear that $A$ contains a subgroup of order $q^3$. Thus, replacing if necessary $A$ by such a subgroup we may assume that $A$ has order $q^3$. In what follows  $A_1,\dots,A_s$ denote the subgroups of index $q$ in $A$.

It was proved in \cite[2.11]{Sh} that under this hypothesis the subgroup $F(G)$ has $(q,m)$-bounded index in $G$.  Hence if $q\leq m$, the subgroup $F(G)$ (and consequently $F_2(G)$) has $m$-bounded index. We will therefore assume that $q>m$. In this case $A$ acts trivially on $C_G(a)/F(C_G(a))$ for any $a\in A^{\#}$. Consequently $[C_G(a),A]\leq F(C_G(a))$ for any $a\in A^{\#}$. 

Observe that $\langle[C_G(A_i),A],[C_G(A_j),A]\rangle$ is nilpotent for any $1\leq i,j\leq s$. This is because the intersection $A_i\cap A_j$ contains a nontrivial element $a$ and the subgroups $[C_G(A_i),A]$ and $[C_G(A_j),A]$ are both contained in the nilpotent subgroup $[C_G(a),A]$.

\begin{lemma}\label{nilp}
 The subgroup $[G,A]$ is nilpotent.
\end{lemma}
\begin{proof}
We argue by contradiction. Suppose $G$ is a counterexample of minimal possible order. By Lemma \ref{solu} the subgroup $[G,A]$ is soluble. Let $V$ be a minimal $A$-invariant normal subgroup of $G$. Then $V$ is an elementary abelian $p$-group and $G/V$ is an $r$-group for some primes $p\neq r$. Write $G=VH$ where $H$ is an $A$-invariant Sylow $r$-subgroup such that $H=[H,A]$. Lemma \ref{L1} says that $H$ is generated by the subgroups $[C_H(A_i),A]$. Thus, $H$ centralizes $[V,A]$ since $[C_V(A_i),A]$ and $[C_H(A_j),A]$ have coprime order for each $1\leq i,j\leq s$. Hence $[V,A]\leq Z(G)$ and by the minimality we conclude that $[V,A]=1$ and $V=C_V(A)$. But then by Lemma \ref{lz} $V\leq Z(G)$ since $V$ is a normal subgroup and $G=[G,A]$. This is a contradiction and the lemma is proved.
\end{proof}

We can now easily complete the proof of Theorem \ref{main2}. By the above lemma $A$ acts trivially on the quotient $G/F(G)$. Therefore $G=F(G)C_G(A)$.
This shows that $F(C_G(A))\leq F_2(G)$. Since the index of $F(C_G(A))$ in $C_G(A)$ is at most $m$, the result follows.

\section{Proof of Theorem \ref{main1}} 

We say that a finite group $G$ is metanilpotent if $\gamma_{\infty}(G)\leq F(G)$.

The following elementary lemma will be useful (for the proof see for example \cite[Lemma 2.4]{AST}).

\begin{lemma}\label{I3} Let $G$ be a metanilpotent finite group. Let $P$ be a Sylow $p$-subgroup of $\gamma_{\infty} (G)$ and $H$ be a Hall $p'$-subgroup of G. Then $P=[P,H]$.
\end{lemma}

Let us now assume the hypothesis of Theorem \ref{main1}. Thus, $A$ is an elementary abelian group of order at least $q^3$ acting on a finite $q'$-group $G$ in such a manner that $\gamma_{\infty} (C_{G}(a))$ has order at most $m$ for any $a\in A^{\#}$. We wish to show that $\gamma_{\infty} (G)$ has $m$-bounded order. Replacing if necessary $A$ by a subgroup we may assume that $A$ has order $q^3$. Since $\gamma_{\infty}(C_{G}(a))$ has order at most $m$, we obtain that $F(C_G(a))$ has index at most $m!$ (see for example \cite[2.4.5]{khukhro}). By \cite[Theorem 1.1]{Eme2} $\gamma_{\infty} (G)$ has $(q,m)$-bounded order. Without loss of generality we will assume that $m!<q$. In particular, $[G,A]$ is nilpotent by Lemma \ref{nilp}.

\begin{lemma}\label{gam}
If $G$ is soluble, then $\gamma_{\infty}(G)=\gamma_{\infty}(C_G(A))$.
\end{lemma}
\begin{proof} We will use induction on the Fitting height $h$ of $G$.

Suppose first that $G$ is metanilpotent. Let $P$ be a Sylow $p$-subgroup of $\gamma_{\infty}(G)$ and $H$ a Hall $A$-invariant $p'$-subgroup of $G$. By Lemma \ref{I3} we have $\gamma_{\infty}(G)=[P,H]=P$. It is sufficient to show that $P\leq\gamma_{\infty}(C_G(A))$. Therefore without loss of generality we assume that $G=PH$. With this in mind, observe that $\gamma_{\infty}(C_G(a))=[C_P(a),C_H(a)]$ for any $a\in A^{\#}$. 

We will prove that $P=[C_P(A),C_H(A)]$. Note that $A$ acts trivially on $\gamma_{\infty}(C_G(a))$ for any $a\in A^{\#}$ since $m<q$. Hence $\gamma_{\infty}(C_G(a))\leq C_P(A)$ for any $a\in A^{\#}$. Let $a,b\in A$. We have $[\gamma_{\infty}(C_G(a)),C_H(b)]\leq[C_P(A),C_H(b)]\leq\gamma_{\infty}(C_G(b))$. Let us show that $P=C_P(A)$.

Assume first that $P$ is abelian. Observe that the subgroup $N=\prod_{a\in A^{\#}}\gamma_{\infty}(C_G(a))$ is normal in $G$. Since $N$ is $A$-invariant, we obtain that $A$ acts on $G/N$ in such a way that $C_G(a)$ is nilpotent for any $a\in A^{\#}$. Thus $G/N$ is nilpotent by \cite{War}. Therefore $P=\prod_{a\in A^{\#}}\gamma_{\infty}(C_G(a))$. In particular $P=C_P(A)$.

Suppose now that $P$ is not abelian. Consider the action of $A$ on $G/\Phi(P)$. By the above, $P/\Phi(P)=C_P(A)\Phi(P)/\Phi(P)$, which implies that $P=C_P(A)$.

Since $P=C_P(A)$ is a normal subgroup of $G$, by Lemma \ref{lz} we deduce that $[H,A]$ centralizes $P$. Therefore, $P=[C_P(A),C_H(A)]$ since $H=[H,A]C_H(A)$. This completes the proof for metanilpotent groups.

If $G$ is soluble and has Fitting height $h>2$, we consider the quotient group $G/\gamma_{\infty} (F_2(G))$ which has Fitting height $h-1$. Clearly $\gamma_{\infty} (F_2(G))\leq \gamma_{\infty}(G)$. Hence, we deduce that $\gamma_{\infty}(G)=\gamma_{\infty}(C_G(A))$.
\end{proof} 

Recall that under our assumptions $[G,A]$ is nilpotent and $C_G(A)$ has a normal nilpotent subgroup of index at most $m!$. Let $R$ be the soluble radical of $G$. Since $G=[G,A]C_G(A)$, the index of $R$ in $G$ is at most $m!$. Lemma \ref{gam} shows that the order of $\gamma_{\infty}(R)$ is at most $m$. We pass to the quotient $G/\gamma_{\infty}(R)$ and without loss of generality assume that $R$ is nilpotent. If $G=R$, we have nothing to prove. Therefore assume that $R<G$ and use induction on the index of $R$ in $G$. Since $[G,A]\leq R$, it follows that each subgroup of $G$ containing $R$ is $A$-invariant. If $T$ is any proper normal subgroup of $G$ containing $R$, by induction the order of $\gamma_{\infty}(T)$ is $m$-bounded and the theorem follows. Hence, we can assume that $G/R$ is a nonabelian simple group. We know that $G/R$ is isomorphic to a quotient of $C_G(A)$ and so, being simple, $G/R$ has order at most $m$.

As usual, given a set of primes $\pi$, we write $O_\pi(U)$ to denote the maximal normal $\pi$-subgroup of a finite group $U$. Let $\pi=\pi(m!)$ be the set of primes at most $m$. Let $N=O_{\pi'}(G)$. Our assumptions imply that $G/N$ is a $\pi$-group and $N\leq F(G)$. Thus,  by the Schur-Zassenhaus theorem \cite[Theorem 6.2.1]{GO} the group $G$ has an $A$-invariant $\pi$-subgroup $K$ such that $G=NK$. Let $K_0=O_\pi(G)$.

Suppose that $K_0=1$. Then $G$ is a semidirect product of $N$ by $K=C_K(A)$. For an automorphism $a\in A^\#$ observe that $[C_N(a),K]\leq\gamma_\infty(C_G(a))$ since $C_N(a)$ and $K$ have coprime order. On the one hand, being a subgroup of $\gamma_\infty(C_G(a))$, the subgroup $[C_N(a),K]$ must be a $\pi$-group. On the one hand, being a subgroup of $N$, the subgroup $[C_N(a),K]$ must be a $\pi'$-group. We conclude that $[C_N(a),K]=1$ for each $a\in A^\#$. Since $N$ is a product of all such centralizers $C_N(a)$, it follows that $[N,K]=1$. Since $K_0=1$ and $K$ is a $\pi$-group, we deduce  that $K=1$ and so $G=N$ is a nilpotent group.

In general $K_0$ does not have to be trivial. However considering the quotient $G/K_0$ and taking into account the above paragraph we deduce that $G=N\times K$. In particular, $\gamma_\infty(G)=\gamma_\infty(K)$ and so without loss of generality we can assume that $G$ is a $\pi$-group. It follows that the number of prime divisors of $|R|$ is $m$-bounded and we can use induction on this number. It will be convenient to prove our theorem first under the additional assumption that $G=G'$.

Suppose that $R$ is an $p$-group for some prime $p\in\pi$. Note that if $s$ is a prime different from $p$ and $H$ is an $A$-invariant Sylow $s$-subgroup of $G$, then in view of Lemma \ref{gam} we have $\gamma_{\infty}(RH)\leq\gamma_{\infty}(C_G(A))$ because $RH$ is soluble. We will require the following observation about finite simple groups (for the proof see for example \cite[Lemma 3.2]{Eme2}).

\begin{lemma}\label{uu}
Let $D$ be a nonabelian finite simple group and $p$ a prime. There exists a prime $s$ different from $p$ such that $D$ is generated by two Sylow $s$-subgroup.
\end{lemma}

In view of Lemma \ref{uu} and the fact that $G/R$ is simple we deduce that $G/R$ is generated by the image of two Sylow $s$-subgroup $H_1$ and $H_2$ where $s$ is a prime different from $p$. Both subgroups $RH_1$ and $RH_2$ are soluble and $A$-invariant since $[G,A]\leq R$. Therefore both $[R,H_1]$ and $[R,H_2]$ are contained in $\gamma_{\infty }(C_G(A))$. 

Let $H=\langle H_1,H_2\rangle$. Thus $G=RH$. Since $G=G'$, it is clear that $G=[R,H]H$ and $[R,G]=[R,H]$. We have $[R,H]=[R,H_1][R,H_2]$ and therefore the order of $[R,H]$ is $m$-bounded. Passing to the quotient $G/[R,G]$ we can assume that $R=Z(G)$. So we are in the situation where $G/Z(G)$ has order at most $m$. By a theorem of Schur the order of $G'$ is $m$-bounded as well (see for example \cite[2.4.1]{khukhro}). Taking into account that $G=G'$ we conclude that the order of $G$ is $m$-bounded.

Suppose now that $\pi(R)=\{p_1,\dots,p_t\}$, where $t\geq2$. For each $i=1,\dots,t$ consider the quotient $G/O_{p_i'}(G)$. The above paragraph shows that the order of $G/O_{p_i'}(G)$ is $m$-bounded. Since also $t$ is $m$-bounded, the result follows.

Thus, in the case where $G=G'$ the theorem is proved. Let us now deal with the case where $G\neq G'$. Let $G^{(l)}$ be the last term of the derived series of $G$. The previous paragraph shows that $|G^{(l)}|$ is $m$-bounded. Consequently, $|\gamma_{\infty}(G)|$ is $m$-bounded since $G/G^{(l)}$ is soluble and $G^{(l)}\leq\gamma_{\infty}(G)$. The proof is now complete.

\end{document}